\documentclass[11pt]{amsart}
\usepackage{verbatim,amssymb,amsmath,mathrsfs,graphicx,hyperref}
\usepackage{xcolor}
  \scrollmode

\newcommand{\R}{{\mathbb R}}
\newcommand{\N}{{\mathbb N}}

\newcommand{\EE}{{\mathbb E}}
\newcommand{\PP}{{\mathbb P}}

\newcommand{\Xh}{\widehat X}
\newcommand{\eul}{{\widehat X^\text{E}}}
\newcommand{\mil}{{\widehat X^\text{M}}}

\newcommand{\eps}{\varepsilon}

\theoremstyle{plain}
\newtheorem{theorem}{Theorem}

\theoremstyle{definition}
\newtheorem{rem}{Remark}

\begin{document}
	
	\title[Complexity of strong approximation of SDE{s} with non-Lipschitz drift coefficient]{On the complexity of strong approximation of stochastic differential equations with a non-Lipschitz drift coefficient}

	\author[M\"uller-Gronbach]
	{Thomas M\"uller-Gronbach}
	\address{
		Faculty of Computer Science and Mathematics\\
		University of Passau\\
		Innstrasse 33 \\
		94032 Passau\\
		Germany} \email{thomas.mueller-gronbach@uni-passau.de}
	
	\author[Yaroslavtseva]
	{Larisa Yaroslavtseva}
	\address{
		621 Institute of Mathematics and Scientific Computing\\
		University of Graz\\
		Heinrichstra{\ss}e 36 \\
		8010 Graz\\
		Austria} \email{larisa.yaroslavtseva@uni-graz.at}

	\begin{abstract} 
		We survey recent developments in the field of
		%discuss the 
		complexity of pathwise approximation in $p$-th mean of the solution of a stochastic differential equation at the final time based on finitely many evaluations of the driving Brownian motion. First, we briefly review the case of equations with globally Lipschitz continuous coefficients, for which an error rate of at least $1/2$ in terms of the number of evaluations of the driving Brownian motion is always guaranteed by using the equidistant Euler-Maruyama scheme. Then we illustrate that giving up the global Lipschitz continuity of the coefficients may lead to a non-polynomial decay of the error  for the Euler-Maruyama scheme or even to an arbitrary slow decay of the smallest possible error that can be achieved on the basis of  finitely many evaluations of the driving Brownian motion. Finally, we turn to recent positive results for equations with a drift coefficient that is not globally Lipschitz continuous. Here we focus on scalar equations with a Lipschitz continuous diffusion coefficient and a drift coefficient that satisfies piecewise smoothness assumptions or has fractional Sobolev regularity and we present corresponding complexity results.

	\end{abstract}

	\maketitle

	\section{Introduction}\label{intro}
	
	We consider a $d$-dimensional system of autonomuous stochastic differential equations (SDEs) on the unit interval,
	\begin{equation}\label{ee1}
		\begin{aligned}
			dX_t &= a(X_t) \, dt + b(X_t) \, dW_t,
			\quad t\in[0, 1],\\
			X_0 &=x_0\in\R^d
		\end{aligned}
	\end{equation}
	with  drift coefficient $a\colon\R^d\to \R^d$, diffusion coefficient $b\colon \R^d\to \R^{d\times m}$ and  $m$-dimensional driving Brownian motion $W$.
	
	The computational goal is to approximate the solution $X$ of the equation~\eqref{ee1}, if it exists, at the final time $t=1$ based on finitely many evaluations of the driving Brownian motion $W$, i.e. by an approximation of the type 
	\begin{equation}\label{ee2}
		\Xh_{n,1} = u(W_{t_1}, \ldots, W_{t_n}),
	\end{equation}
	where $u\colon \R^{m\times n}\to\R^d$ is a measurable function and $t_1,\dots,t_n\in [0,1]$. The quality of such an approximation is measured by its $L_p$-error, i.e.
	\[
	\EE \bigl[\|X_1-\widehat X_{n,1}\|^p\bigr]^{1/p}.
	\]
	
		There are two types of results we want to discuss. First, we are interested in upper error bounds for specific methods $\Xh_{n,1}$ of the type~\eqref{ee2} in terms of the number $n$ of evaluations of $W$ that are used by $\Xh_{n,1}$, e.g.  
	\begin{equation}\label{ee3}
		\begin{aligned}
			\EE \bigl[\|X_1-\widehat X_{n,1}\|^p\bigr]^{1/p}\preceq n^{-\alpha}
		\end{aligned}
	\end{equation}
	with $\alpha\in (0,\infty)$, where for two sequences $(a_n)_{n\in\N} $ and $(b_n)_{n\in\N} $ of nonnegative reals we use the notation $a_n\preceq b_n$ if there exists $c\in (0,\infty)$ such that $a_n\le c b_n$  for all $n\in\N$. Furthermore, we write $a_n\succeq b_n$ if $b_n\preceq a_n$ and we write $a_n\asymp b_n$ if $a_n\preceq b_n\preceq a_n$. 
	
	Second, we are interested in lower error bounds that hold uniformly for every approximation of the type~\eqref{ee2}, e.g.
	\begin{equation}\label{ee4}
		e_n^{(p)} := \inf_{
			\substack{
				t_1,\dots,t_n\in [0,1],
				\\
				u\colon \R^{m\times n}\to \R^d \text{ measurable}
			}
		}
		\EE\bigl[ \|X_1-u(W_{t_1}, \ldots,
		W_{t_n})\|^p\bigr]^{1/p}  \succeq 
		n^{-\beta}
	\end{equation}
	with $\beta\in (0,\infty)$. In \eqref{ee4}, the quantity $e_n^{(p)} $ is the smallest possible $L_p$-error  that can be achieved based on $n$ evaluations of $W$ at fixed time points $t_1,\dots,t_n\in [0,1]$ and is called the  $n$-th minimal error in $p$-th mean sense for non-adaptive methods. 
	
	Ideally, the bounds \eqref{ee2} and \eqref{ee3} are sharp in the sense of $\alpha =\beta$. In this case we have
	\[ 
	\EE \bigl[\|X_1-\widehat X_{n,1}\|^p\bigr]^{1/p} \asymp e_n^{(p)}\asymp n^{-\beta},
	\]
	i.e. $\beta$ is the optimal $L_p$-error rate and the sequence of methods $(\widehat X_{n,1})_{n\in\N}$ performs asymptotically optimal.
	
	%T could be carried out formally. I suggest not to do that.
	More general, one can consider the class of all adaptive methods that are based on $n$ sequential evaluations of $W$ on average and study the asymptotic behaviour of the corresponding $n$-th minimal errors in $p$-th mean. See e.g.~\cite{hhmg2019,m04} for a formal definition of this class of methods. In the present survey, we will mention the available corresponding results only in an informal way.
	
		There are two particular methods that play a prominent role in the literature on numerical approximation of SDEs, namely the Euler-Maruyama scheme, introduced in 1955 in~\cite{m55}, and the Milstein scheme, introduced in 1974 in~\cite{m74}. For $n\in \N$, the Euler-Maruyama scheme with $n$ equidistant time steps is given by $ \eul_{n,0} =x_0$ and
	\begin{equation}
	\begin{aligned}
	\eul_{n, (\ell+1)/n}&=
	\eul_{n,\ell/n}
	+
	a(  \eul_{n,\ell/n})
	\, \frac{1}{ n }
	+
	b(\eul_{n,\ell/n}
	)\, (W_ {(\ell+1)/n}-W_ {\ell/n})
\end{aligned}
\end{equation}	
for $\ell = 0,\dots,n-1$. Thus, $\eul_n$  generalizes the Euler scheme for ordinary differential equations by adding the simple approximation 	$b(\eul_{n,\ell/n}
)(W_ {(\ell+1)/n}-W_ {\ell/n})$ of  the stochastic integral $\int_{t_\ell}^{t_{\ell+1}} b(X_s)\, dW_s$ in each step.

The Milstein scheme with $n$ equidistant time steps requires differentiability of the diffusion coefficient $b$ and employs a more accurate approximation of the stochastic integrals $\int_{t_\ell}^{t_{\ell+1}} b(X_s)\, dW_s$. It is given by $ \mil_{n,0} =x_0$ and
	\begin{equation}
	\begin{aligned}
 				\mil_{n, (\ell+1)/n}&=
				\mil_{n, \ell/n}
				+
				a(  \mil_{n,\ell/n})
				\, \frac{1}{ n }
				+
				b(\mil_{n,\ell/n}
				)\, (W_ {(\ell+1)/n}-W_ {\ell/n})\\
				&\qquad\qquad +\sum_{i, j=1}^m \nabla b^{(i)}b^{(j)}(\mil_{n, \ell/n})\, \int_{\ell/n}^{(\ell+1)/n} (W^{(i)}_{ s}-W^{(i)}_{ \ell/n})\, dW^{(j)}_{ s}			
			\end{aligned}
		\end{equation}
for $\ell = 0,\dots,n-1$, where $b^{(i)}=(b_{1,i}, \ldots, b_{d,i})^T$ denotes the $i$-th column of the matrix-valued mapping $b$ and $W^{(i)}$ denotes the $i$-th component of $W$ for $i=1,\dots,m$. 

\begin{rem}\label{rem1}
Note that the Milstein approximation $\mil_{n, 1}$ of $X_1$ uses not only the evaluations $W_{1/n},\dots,W_1$ of the driving Brownian motion but also the iterated It\^{o} integrals 
\[
J^{(\ell)}_{i,j} =  \int_{\ell/n}^{(\ell+1)/n} (W^{(i)}_{ s}-W^{(i)}_{\ell/n})\, dW^{(j)}_{ s}
\]
for $i,j\in\{1,\dots,m\}$ and $\ell\in\{0,\dots,n-1\}$ and is therefore not a method of type~\eqref{ee2} in general. 

However, note that
		\[
J^{(\ell)}_{i,j} + J^{(\ell)}_{j,i}	  = \begin{cases} W^{(i)}_{ (\ell+1)/n}-W^{(i)}_{\ell/n})\, (W^{(j)}_{(\ell+1)/n}-W^{(j)}_{\ell/n}), & \text{if }i\neq j,\\ (W^{(i)}_{ (\ell+1)/n}-W^{(i)}_{\ell/n})^2 -1/n, & \text{if }i = j \end{cases}
\]
for all $i,j\in \{1,\dots,m\}$ and all $\ell\in \{0,\dots,n-1\}$. As a consequence we obtain that if the diffusion coefficient $b$ satisfies the commutativity condition 
	\[
		\forall i, j\in\{1, \ldots, m\}\colon \,\, \nabla b^{(i)}b^{(j)}=\nabla b^{(j)}b^{(i)} \tag{C}
	\]
	then, for  all $\ell\in \{0,\dots,n-1\}$,
\[
\sum_{i, j=1}^m \nabla b^{(i)}b^{(j)}(\mil_{n, \ell/n})\, \int_{\ell/n}^{(\ell+1)/n} (W^{(i)}_{ s}-W^{(i)}_{ \ell/n})\, dW^{(j)}_{ s} = g_{\ell+1}(W_{1/n}, \ldots, W_{(\ell+1)/n})
\]
with a suitable function $g_{\ell+1}\colon \R^{(\ell+1)\times n} \to \R$. Thus, if $b$ satisfies (C) then $\mil_{n, 1}$ is a method of type~\eqref{ee2}. Note that the condition (C) is in particular satisfied if $m=1$, i.e. if the driving Brownian motion $W$ is one-dimensional.  
\end{rem}  

%%%%%%%%%%%%%%%%%%%%%

	\section{Complexity results for Lipschitz coefficients}\label{Lip}
	
 We first present the classical  upper error bound for the Euler-Maruyama scheme in the case of globally Lipschitz continuous coefficients, see e.g. ~\cite{kp92} for a proof.  
	
	\begin{theorem}\label{t1}  
		Assume that $a$ and $b$ are Lipschitz continuous. Then, for all  $p\in [1,\infty)$,
	\[
	\EE \bigl[\|X_1-\eul_{n, 1}\|^p\bigr]^{1/p} \preceq n^{-1/2}.
	\]
	\end{theorem}

If the coefficients $a$ and $b$ satisfy stronger smoothness assumptions than Lipschitz continuity we may even achieve an $L_p$-error rate of at least $1$ by using the Milstein approximation. 
	\begin{theorem}\label{t2}  	
	Assume that $a$ and $b$ are differentiable with bounded, Lipschitz continuous derivatives. Then, for all  $p\in [1,\infty) $,
\[
		\EE \bigl[\|X_1-\mil_{n,1}\|^p\bigr]^{1/p} \preceq n^{-1}.
\]
		\end{theorem}
	
See~\cite{MG02_habil} for a proof of Theorem~\ref{t2} and~\cite{m74, Faure92, hmr01} for proofs of previous versions of Theorem~\ref{t2} using more restrictive assumptions on the coefficients $a$ and $b$. 

In case of a scalar driving Brownian motion, the Milstein approximation $\mil_{n,1}$ is of type~\eqref{ee2}, see Remark~\ref{rem1}. Hence Theorem~\ref{t2} implies that in this case  an $L_p$-error rate of at least $1$ can be achieved based on finitely many evaluations of the driving Brownian motion. The following result, established in \cite{hhmg2019}, shows that under very mild local conditions on the coefficients $a$ and $b$ the error rate $1$ can not be improved in general.

\begin{theorem}\label{t3}
Let $d=m=1$ and assume that equation~\eqref{ee1} has a strong solution  $X$. Let $ I\subset\R$ be open,  let $t_0\in [0,1]$ and assume that
	\begin{align*}
		\text{(i)}& \quad a\text{ and }b \text{ are three times continuously differentiable on }I,\hspace*{17cm}\\[.05cm]
		\text{(ii)}&\quad \forall\, x\in I\colon\,\, b(x)\neq 0 \wedge (a'b-ab' -\tfrac{1}{2} b^2 b'')(x)\neq 0,\\[.07cm]
		\text{(iii)}&\quad \PP(X_{t_0}\in I) > 0.
	\end{align*}
Then 
	\begin{equation}\label{c1}
	e_n^{(1)} \succeq n^{-1}.
	\end{equation}
\end{theorem} 

\begin{rem}\label{rem2}
		The second part of condition (ii) in Theorem~\ref{t3} is motivated by~\cite[Theorem 2.1]{Yamato}, which 
	shows that if the Lie algebra associated with the coefficients of the equation~\eqref{ee1} in Stratonovich form is nilpotent then the solution of the equation~\eqref{ee1} has a finite It\^{o}-Taylor expansion. In particular, this result has the following consequence:  if $d=m=1$, 
 $a,b\in C^\infty(\R;\R)$ and $a'b-ab' -\tfrac{1}{2} b^2 b'' = 0$, i.e. the Lie algebra associated with the Stratonovich coefficients $a-\tfrac{1}{2}b'b$ and $b$ is nilpotent of step $1$, then there exists a measurable function $f\colon \R \to \R$ such that a.s.
		\[
		X_1= f(W_1).
		\]
In this case, $X_1$ can be approximated with error zero by using one evaluation of the driving Brownian motion $W$ and thus the lower bound \eqref{c1} can not be true. A simple example for this situation is provided by the case of a geometric Brownian motion, i.e. 
\[
a(x) = \alpha x,\, b(x) = \beta x
\]
for all $x\in \R$ with $\alpha,\beta\in\R$. Here $(a'b-ab' -\tfrac{1}{2}b^2b'')(x) = \alpha\beta x - \alpha x \beta =0$ for all $x\in \R$ and 
\[
X_1= \exp(\alpha-\beta^2/2 + \beta W_1).
\]
\end{rem}  

In the case when the diffusion coefficient $b$ does not satisfy the commutativity condition (C), the Milstein approximation $\mil_{n,1}$ need not to be of type~\eqref{ee2} and therefore Theorem~\ref{t2} can not be used to obtain an $L_p$-error rate of at least $1$ based on finitely many evaluations of $W$ at points in $[0,1]$. In fact, the following result shows that in this case the $L_p$-error rate $1/2$ achieved by the Euler-Maruyama scheme, see Theorem~\ref{t1}, can not be improved for $p\in [2,\infty)$.

\begin{theorem}\label{t4}
	Assume that $a$ and $b$ have bounded, Lipschitz continuous derivatives and that the commutativity condition (C)  does not hold a.s. on the trajectories of $X$, i.e.  
	\begin{equation}\label{noncom}
	%\EE\Bigl [ \int_0^1\sum_{1\le i < j\le m} \bigl \| \bigl( \nabla b^{(i)}b^{(j)}  - \nabla %b^{(j)}b^{(i)}\bigr)(X_t)\bigr\|\, dt\Bigr] >0.
	\PP\Bigl(\Bigl \{ \forall t\in [0,1]\colon \sum_{1\le i < j\le m} \bigl \| \bigl( \nabla b^{(i)}b^{(j)}  - \nabla b^{(j)}b^{(i)}\bigr)(X_t)\bigr\| =0\, \Bigr\}\Bigr) <1.
	\end{equation}
	Then 
		\[
	e_n^{(2)} \succeq n^{-1/2}.
	\]
\end{theorem}
See~\cite[Chapter IV.2]{MG02_habil} for a proof of Theorem~\ref{t4}. A simple example for the setting of Theorem~\ref{t4} is provided by the $2$-dimensional equation 
\[
		dX^{(1)}_{t} = dW^{(1)}_{t},\,\, dX^{(2)}_{t} = X^{(1)}_{t}\, dW^{(2)}_{t}
\]
with a $2$-dimensional driving Brownian motion $W=(W^{(1)},W^{(2)})$ and initial value $x_0 = 0\in \R^2$. Thus, 
\[
X^{(1)}_{t} = W^{(1)}_{t},\, \,X^{(2)}_{t} = \int_0^t W^{(1)}_{s}\, dW^{(2)}_{s}
\]
for all $t\in [0,1]$. Here, $a=0$ and the diffusion coefficient $b\colon \R^2\to \R^{2\times 2},\, (x_1,x_2)^\top \mapsto \begin{pmatrix} 1 & 0\\ 0 & x_1\end{pmatrix} $ satisfies 
\[
\nabla b^{(1)}b^{(2)}  - \nabla b^{(2)}b^{(1)} = (0,-1)^\top.
\]
By Theorem~\ref{t4}  we conclude that
\[
\inf_{t_1,\dots,t_n\in [0,1]} \EE\Bigl[\Bigl| \int_0^1 W^{(1)}_{t}\, dW^{(2)}_t\, dt - \EE\Bigl[ \int_0^t W^{(1)}_{t}\, dW^{(2)}_{t}\,dt\, \bigl| W_{t_1},\dots, W_{t_n}\Bigr] \Bigr|^2\Bigr] \succeq n^{-1},
\]
which was previously known from~\cite{ClarkCameron1980}.

Combining Theorem~\ref{t1} with Theorem~\ref{t4} and Theorem~\ref{t2} with Theorem~\ref{t3} we obtain the following complexity result.

\begin{theorem}\label{t5} Assume that $a$ and $b$ have bounded, Lipschitz continuous derivatives.
\begin{itemize}
	\item[(i)] If $b$ satisfies~\eqref{noncom}, i.e. the commutativity condition (C)  does not hold a.s. on the trajectories of $X$, then 
	\[
	\forall p\in[2,\infty)\colon\,\, e_n^{(p)} \asymp n^{-1/2}.
	\]
		\item[(ii)] If $b$ satisfies the commutativity condition (C) then 
		\[
	\forall p\in[1,\infty)\colon\,\, e_n^{(p)} \preceq n^{-1}.
	\]
	\item[(iii)] If $d=m=1$ and (i) to (iii) in Theorem~\ref{t3} are satisfied then 
	\[
	\forall p\in[1,\infty)\colon\,\, e_n^{(p)} \asymp n^{-1}.
	\]
\end{itemize}
\end{theorem}

\begin{rem}\label{rem3}

\noindent	
\begin{itemize}
	\item[(i)] The lower bounds in Theorem~\ref{t3} and Theorem~\ref{t4} are  also valid for the class of all adaptive methods based on {\small $n$} sequential evaluations of {\small $W$} on average, see~\cite{hhmg2019} and~\cite{MG02_habil}, respectively. Thus, in the setting of these two theorems, the  respective $L_p$-error rates  $1$ and $1/2$ can not be improved by using e.g. an adaptive step-size control. This situation can change drastically if the Lipschitz continuity of the drift coefficient is abandoned, see Remark~\ref{rem7}. 
	\item[(ii)] Choosing an adaptive step-size control can, however, lead to much smaller asymptotic constants in the asymptotic behaviour of the resulting $L_p$-errors compared to what can be achieved by evaluating $W$ at equidistant time points. See~\cite{MG02_habil}  and~\cite{m04} for a detailed analysis of this issue. 
\end{itemize}
\end{rem}

\section{Negative results for non-Lipschitz coefficients }\label{nonLipneg}	
	
The first negative result for $L_p$-approximation of SDEs  with coefficients that are not Lipschitz continuous seems to have been obtained in~\cite{hjk11}, where  equations with coefficients of superlinear growth are studied, which are of  importance in many fields of applications. It turns out that in this situation not even convergence in the mean can be guaranteed for the Euler-Maruayama scheme with equidistant steps. Even worse, the following result from~\cite{hjk11} essentially shows that if either the drift coefficient grows in a higher polynomial order than linearly and the diffusion coefficient grows slower than that, or the diffusion coefficient grows in a higher polynomial order than linearly and the drift coefficient grows slower than that then the mean error of the equidistant Euler-Maruyama scheme  tends to infinity with growing number of steps.

\begin{theorem}\label{t6}
Let $d=m=1$ and assume that the equation~\eqref{ee1} has a solution $X$ with $\EE[|X_1|]<\infty$. If $b(x_0)\neq 0$ and there exist $C\in (0,\infty)$, $\alpha\in [0,\infty)$ and $\beta\in (\max(1,\alpha),\infty)$ such that either
\[
\forall x\in (-\infty,C] \cup [C,\infty)\colon\,\, |a(x)| \ge |x|^\beta/C\text{ and }|b(x)|\le C|x|^\alpha 
\]
or
\[
\forall x\in (-\infty,C] \cup [C,\infty)\colon\,\, |b(x)| \ge |x|^\beta/C\text{ and }|a(x)|\le C|x|^\alpha 
\]
then
\[
\lim_{n\to\infty}
\EE \bigl[|X_1-\eul_{n,1}|\bigr] =\infty.
\]
\end{theorem}

A simple example for the situation of Theorem~\ref{t6} is given by the equation $dX_t = -X_t^3\,dt + dW_t$, which satisfies the first condition in Theorem~\ref{t6} with the choice $C=1$, $\alpha=0$ and $\beta=3$. A practically more relevant example is provided by the Ginzburg-Landau equation given by 
\[
dX_t = \bigl( (\eta +\gamma^2/2)X_t -\lambda X_t^3\bigr)\, dt + \gamma X_t \, dW_t
\]
with initial value $x_0\in (0,\infty)$ and parameters $\eta\in [0,\infty)$ and $\gamma,\lambda\in (0,\infty)$, which is used to model phase transitions in the theory of superconductivity. This equation  satisfies the first condition in Theorem~\ref{t6} with the choice $\alpha=1$, $\beta= 3$ and $C=\max(1,\gamma, 2/\lambda, (2\eta+\gamma^2)/\lambda)$, see~\cite{hjk11}.

The following result from~\cite{hhj12} shows that a bad performance in $p$-th mean of the Euler-Maruyama scheme may also happen if the coefficients of the equation~\eqref{ee1} are bounded and infinitely often differentiable.

\begin{theorem}\label{t7}
Let $d=4$ and $m=1$. There exist a bounded $a\in C^{\infty}(\R^4,\R^4)$ and a  bounded, Lipschitz continuous $b\in C^{\infty}(\R^4,\R^4)$ such that 
	\[
	\forall\alpha \in (0,\infty)\colon\,\,	\lim_{n\to\infty} n^\alpha\cdot 
	\EE \bigl[\|X_1-\eul_{n,1}\|\bigr] =\infty.	
	\] 
\end{theorem}
In the setting of Theorem~\ref{t7}, the equidistant Euler-Maruyama scheme does not converge in the mean with no polynomial decay. The next result from~\cite{JMGY15} shows that this bad performance is not a specific problem of only the equidistant Euler-Maruyama scheme but may even be the case for any approximation of the type~\eqref{ee2}. Even worse, the best possible speed of convergence in the mean that can be achieved by approximations of type~\eqref{ee2} can be arbitrarily slow.

\begin{theorem}\label{t8}
	Let $d=4$ and $m=1$. 
For every  sequence $(\varepsilon_n)_{n\in\N}$ in $(0,\infty)$ with $\lim_{n\to\infty} \varepsilon_n = 0$ 	
		there exist a bounded $a\in C^{\infty}(\R^4,\R^4)$ and a  bounded, Lipschitz continuous $b\in C^{\infty}(\R^4,\R^4)$ such that
			\[
			\qquad \forall n\in\N\colon\,\, e^{(1)}_n \geq 
					\varepsilon_n.
			\]
			\end{theorem}

\begin{rem}\label{rem3b}
	\noindent
	\begin{itemize}
		\item[(i)] The statement in Theorem~\ref{t8} extends to the class of all adaptive methods based on $n$ sequential evalutions of $W$ on average, see~\cite{Y17}. 
		\item[(ii)] In~\cite{GJS17} it is shown that the negative result in Theorem~\ref{t8} still holds true in dimension $d=2$. Whether the result is also valid in dimension $d=1$ is so far an open question.
	\end{itemize}
\end{rem}	

A crucial point in the proof of Theorem~\ref{t8} in~\cite{JMGY15} is to construct a drift coefficient $a$ which has a rapidly oscillating component. One might therefore think that imposing a polynomial growth condition on the derivatives of the coefficients $a$ and $b$ could guarantee that a polynomial decay of the error  in $p$-th mean can be achieved. However, the following result from~\cite{MGY2018} shows that in general this believe is wrong, even if the solution of the equation~\eqref{ee1} is required to have a finite uniform norm in the mean. 

\begin{theorem}\label{t9}
Let $d=7$ and $ m=1$. There exist 
$a,b\in C^\infty(\R^7,\R^7)$ such that $b$ is Lipschitz continuous, $a'$ has at most  polynomial growth, 
$ \EE[\|X\|_\infty]<\infty$ and 
\[
	e^{(1)}_n \succeq \frac{1}{\ln^2 (n+1)}.
\]
\end{theorem}
	
Finally, we illustrate, by a result from~\cite{HefJ17}, that an arbitrarily slow polynomial decay of the mean error of any method based on equidistant evaluations of $W$ is possible if the drift coefficient $a$ is Lipschitz continuous but the diffusion coeffient $b$ is only H\"older continous of order $1/2$. 
\begin{theorem}\label{t10}
Let $d=m=1$. Consider the Cox-Ingersoll-Ross (CIR) process given by
		\begin{align*}
			dX_t &=  (\kappa-\lambda X_t)\,dt + \sigma\sqrt{X_t} \, dW_t,
			\quad t\in[0, 1],\\
			%\vspace{-0.1cm}
			X_0&=x_0
		\end{align*}
	with $\kappa, \sigma,\lambda, x_0\in (0,\infty)$. If the Feller index $\nu = 2\kappa/\sigma^2$ satisfies $\nu\in (0,1)$ then 
		\[ 
		\inf_{			u\colon\R^n\to \R \text{ measurable}
		}
		\EE \bigl[|X_1-u(W_{1/n}, W_{2/n}, \ldots,
		W_{1})|\bigr] \succeq n^{-\nu}.
		\]
\end{theorem}	

\section{Positive results for non-Lipschitz coefficients}\label{nonLippos}

Motivated by the negative result Theorem~\ref{t6} on the potential divergence of the Euler-Maruyama scheme in case of SDEs with coefficients of superlinear growth, research on $L_p$-approximation of SDEs with non-Lipschitz coefficients has been significantly intensified since then. Mainly, the following three settings are studied. 
\begin{itemize}
	\item {\bf Locally Lipschitz continuous coefficients} \\[.02cm] 
	Here, new types of methods have been constructed by combining the classical Euler-Maruyama scheme and the Milstein scheme with strategies of taming, projecting, truncating and adaptive step-size control in order to cope with superlinearly  growing coefficients. Typically, the coefficients $a$ and $b$ of the equation~\eqref{ee1} are required to satisfy at least a monotonicity-type condition and a coercivity condition and to be locally Lipschitz continuous with a local Lipschitz constant of at most polynomial growth. Under these conditions, the classical $L_p$-error rates of at least $1/2$ and $1$ are proven for the above mentioned variants of the Euler-Maruyama scheme  and the Milstein scheme, respectively. See~\cite{GanHeWang2020,Hatz2020, hjk12,KuSab2019,Sabanis2013, Sabanis2016, SabZa2019,TretZa2013,WangGan2013,MaZhang2017} for results on taming, \cite{Beynetal2015, Beynetal2017} for results on projecting, \cite{HuLiMao2018,HutzenthalerJentzen2020,Mao2016,YangWuKloedenMao2020} for results on truncating and \cite{FangGiles} for the approach of adaptive step size control. 
	We add that in this setting an  $L_p$-error rate of at least $1/2$ and $1$ is also achieved by implicit versions of the Euler-Maruyama scheme and the Milstein scheme, respectively, see e.g.~\cite{hms02, Hu1996,MaoSzpruch2013,TretZa2013}.\\[-.4cm]
	\item {\bf H\"older continuous coefficients}\\[.02cm]
	Since about 2010, $L_p$-error rates are also obtained for strong approximation of SDEs with H\"older continuous coefficients, see e.g.~\cite{BaHuYu2019,BerkaouiBossyDiop2008, BosOli2018,butdage2021,gerlamlin2023,GyoengyRasonyi2011,MeTa2017,NgoTag2016,    NgoTag2017b}. Roughly speaking, if the drift coefficient $a$ is H\"older continuous of order $\alpha\in (0,1)$ and the diffusion coefficient $b$ is Lipschitz continuous and satisfies further regularity assumptions then an $L_p$-error rate of at least $1/2-$ and $(\alpha+1)/2-$ 
	  is achieved by the Euler-Maruyama scheme and the Milstein scheme, respectively, see~\cite{butdage2021,gerlamlin2023}. For the limiting case of $\alpha=0$, i.e. a measurable and bounded drift coefficient, an $L_p$-error rate of at least $1/2-$ is established in~\cite{DG2020,DGL22}. On the other hand, if the drift coefficient $a$ is Lipschitz continuous (and satisfies further regularity conditions) and the diffusion coefficient is of the form $b(x) =|x|^\beta$ with $\beta\in(1/2,1) $ then an $L_p$-error rate of at least $1/2$ and $1$ is achieved by symmetrized variants of the Euler-Maruyama scheme and the Milstein scheme, respectively, see~\cite{BerkaouiBossyDiop2008, BosOli2018,NeuenkirchMickel2024}. In the  case of the CIR-process, i.e. $a(x) = \kappa - \lambda x$ and $b(x)=\sigma | x|^{1/2}$, see Theorem~\ref{t10}, 
	  %    which serves as a model of state dependent volatility in mathematical finance, 
	  $L_p$-error rates are established for implicit versions of the Euler-Maruyama scheme in~\cite{Alfonsi2013,DereichNeuenkirchSzpruch2012,HeHe17,HutzenthalerJentzenNoll2014CIR,NeuenkirchSzpruch2014}. In~\cite{HH2018}, a  square root truncated Milstein scheme is introduced that achieves an $L_p$-error rate of at least $\min(1/2, 2\kappa/\sigma^2)/p$, which seems to be the best $L_p$-error rate for the CIR process known up to now and is optimal in the class of all methods based on equidistant evaluation of the driving Brownian motion if $2\kappa/\sigma^2<1/2$ and $p=1$, see  Theorem~\ref{t10}. \\[-.4cm]
	  \item {\bf Drift coefficients with piecewise regularity or fractional Sobolev regularity}\\[.02cm]
	  Recently,  polynomial $L_p$-error rates for strong approximation of SDEs are also established in settings that allow for discontinuities of the drift coefficient. For settings with a drift coefficient that is at least piecewise Lipschitz continuous, such rates are proven for variants of the Euler-Maruyama scheme and the Milstein scheme in~\cite{DoNgoPho2024,GanHu2022,LS2016, LS2017, LS2018, MGSY2022, MGY2020, MGY2022, NSS2019,  NgoTag2016,  NgoTag2017b,  NgoTag2017a, SS2022}. See also~\cite{GLN2019} for related extensive  simulation studies.  In case that the drift coefficient has fractional Sobolev regularity of order $s\in (0,1)$ and the diffusion coefficient is sufficiently regular, an $L_p$-error rate of at least $(1+s)/2-$ is proven in~\cite{DGL22,NS2021}
\end{itemize}

We add that  polynomial $L_p$-error rates for strong approximation of SDEs with non-Lipschitz coefficients are meanwhile also obtained in the case of jump-diffusions, see~\cite{PSS2023a,PSS2023b,PS2021},  in the case of a fractional Brownian noise, see~\cite{butdage2021} and in a general L\'evy-noise setting, see~\cite{butdage2024,Kumar2021a}.

We turn to a more detailed presentation of complexity results for the case of a drift coefficient with piecewise regularity and the case of a drift coefficient with fractional Sobolev regularity. Throughout the following, we restrict to the case of scalar SDEs, i.e.  $d=m=1$. 

\subsection{Complexity results for scalar SDEs with a piecewise regular drift coefficient}\label{piecewise}

The following result on the performance of the  Euler-Maruyama scheme extends Theorem~\ref{t1} to the case that the drift coefficient is piecewise Lipschitz continuous.  It shows that in this case, the Euler-Maruyama scheme still achieves, for all $p\in [1,\infty)$,  an $L_p$-error rate of at least $1/2$. See~\cite{MGY2020} for a proof.

	\begin{theorem}\label{t11}  
	Let $d=m=1$, let $k\in\N_0$ and assume that $a$ and $b$ satisfy
	\begin{align}
	\text{There exist }  -\infty = \xi_0 < \xi_1 < \dots < \xi_{k} < \xi_{k+1} = \infty  \text{ such that} \hspace*{.45cm}& \tag{A1}\\
	a \text{ is Lipschitz on } (\xi_{i-1}, \xi_i)  \text{ for }i=1,\dots,k+1,  \hspace*{2.7cm}&\notag\\[.1cm]
	b \text{ is Lipschitz on }\R \text{ and }b(\xi_i)\neq 0 \text{ for }i=1,\dots,k.\hspace*{2.2cm}& \tag{A2}
\end{align}
Then for all  $p\in [1,\infty)$,
	\[
	\EE \bigl[|X_1-\eul_{n, 1}|^p\bigr]^{1/p} \preceq n^{-1/2}.
	\]
\end{theorem}

A simple example for the setting of Theorem~\ref{t11} is given by the equation
\[
dX_t = 1_{[0,\infty)} (X_t)\, dt + dW_t
\]
with initial value $x_0=0$. In this case,
\[
X_1 = \int_0^1 1_{[0,\infty)} (X_t)\, dt + W_1,
\]
i.e. we essentially want to approximate $ \int_0^1 1_{[0,\infty)} (X_t)\, dt$, the total amount of times $t\in[0,1]$ that the solution $X$ stays on the positive side of the reals. 

\begin{rem}\label{rem5}
\noindent	
\begin{itemize}
	\item[(i)] Under the assumptions (A1) and (A2), the SDE~\eqref{ee1} has a unique strong solution, see~\cite{LS2016}.
	\item[(ii)] In the setting of Theorem~\ref{t11}, an $L_p$-error rate of at least $1/2$ is also achieved by a transformed Euler-Maruyama scheme, see~\cite{LS2016}, and an Euler-Maruyama  scheme with adaptive step-size control, see~\cite{NSS2019}. 
	
	In the first paper, a suitable bi-Lipschitz mapping $G\colon\R\to\R$ is constructed, such that the process $Y=(G(X_t))_{t\in[0,1]}$ is the strong solution of an SDE with globally Lipschitz continuous coefficients. The corresponding Euler-Maruyama approximation $\widehat Y_{n,1}$ of $Y_1$ thus satisfies $\EE[|Y_1-\widehat Y_{n,1}|^p] \preceq n^{-p/2}$ by Theorem~\ref{t1}. Moreover, by the Lipschitz continuity of $G^{-1}$, 
	\[
	|X_1-G^{-1}(\widehat Y_{n,1})| =	|G^{-1}(Y_1)-G^{-1}(\widehat Y_{n,1})|  \preceq |Y_1-\widehat Y_{n,1}|.
	\]
	Combining the latter  two estimates yields the $L_p$-error rate $1/2$ for the approximation $G^{-1}(\widehat Y_{n,1})$ of $X_1$.
	
	In the second paper, the step-size of the  Euler-Maruyama scheme is no longer equidistant but adapted to the distance of the actual value of the scheme to the nearest discontinuity of the drift coefficient $a$.   
	\item[(iii)] We briefly illustrate the main difficulty in the proof of the error bound for the Euler-Maruyama approximation under the assumptions in Theorem~\ref{t11} in comparison to the proof of this error bound for global Lipschitz continuous coefficients in Theorem~\ref{t1}. 
	
	For convenience we assume that $b=1$ and that $a$ has a single discontinuity at $\xi=0$ and we restrict to the analysis of the $L_1$-error. Clearly, for every $t\in [0,1]$,
			\[
		X_t = x_0 +\int_0^t a(X_s)\, ds +W_t,\,\, \widehat X^\text{E}_{n,t} = x_0 +\int_0^t a(\widehat X^\text{E}_{n,\underline s_n})\, ds +W_t,  
		\]
		where $(\eul_{n,t})_{t\in[0,1]}$ denotes the time-continuous version of the Euler-Maruyama scheme with $n$ equidistant steps and $\underline s_n = \lfloor sn\rfloor /n$ is the grid point that is closest to $s$ from below.
	Hence,
		\begin{align*}
			X_t - \widehat X^\text{E}_{n,t} &= \int_0^t \bigl(a(X_s) - a(\widehat X^{\text{E}}_{n,\underline s_n})\bigr)\, ds \\
			&= \int_0^t \bigl(a(X_s) - a(\widehat X^{\text{E}}_{n,s})\bigr)\, ds +  \int_0^t \bigl(a(\widehat X^{\text{E}}_{n,s}) - a(\widehat X^{\text{E}}_{n,\underline s_n})\bigr)\, ds.
		\end{align*}
If the drift coefficient $a$ would be Lipschitz continuous we could conclude that
	\[
	\EE\bigl [|X_t - \widehat X^\text{E}_{n,t}|\bigr]  \le c  \int_0^t \EE\bigl [| X_s - \widehat X^{\text{E}}_{n,s}|\bigr]\, ds + d \int_0^t \EE\bigl [|\widehat X^{\text{E}}_{n,s} - \widehat X^{\text{E}}_{n,\underline s_n}|\bigr ] \, ds,
	\]
	where $c\in (0,\infty)$ does not depend on $n$ or on $t$.
Since the time-continuous Euler-Maruyama scheme is H\"older continuous  of order $1/2$ in the mean, we obtain 
	\[
\EE[| X_t - \eul_{n,t}|] \le c  \int_0^t \EE\bigl [| X_s - \widehat X^{\text{E}}_{n,s}|\bigr]\, ds + c n^{-1/2},
\]
	where $c\in (0,\infty)$ does not depend on $n$ or on $t$, and applying  Gronwall's inequality would yield the desired bound.

In the presence of a discontinuity of $a$ at zero, one can use the transformation $G$ mentioned in part (ii) as well as Gronwall's inequality to derive
	\[
\EE\bigl [| X_1 - \eul_{n,1}|\bigr] \preceq n^{-1/2} + D_n,
	\]
where the term
	\[
	D_n = 
	\EE\Bigl[\Bigl | \int_0^1  1_{\{\eul_{n,t}\cdot \eul_{n,\underline{t}_n} \le 0\}}\, dt\Bigr|\Bigr]
	\]
is the average amount of times $t\in[0,1]$ such that the values of the Euler-Maruyama scheme at times $t$ and $\underline t_n$ are on different sides of the point zero, which quantifies the impact of the discontinuity of $a$. Using a  time mirroring technique and occupation time estimates for the time-continuous Euler-Maruyama scheme one can show that $D_n\preceq n^{-1/2}$, see~\cite[Proposition 1]{MGY2020}, which completes the proof of Theorem~\ref{t11}.
\end{itemize}
\end{rem}

Similar to the case of globally Lipschitz continuous coefficients, an  $L_p$-error rate  larger than $1/2$ can be achieved by imposing further piecewise regularity conditions on the coefficients $a$ and $b$. See~\cite{MGY2022} for a proof of the following result. 

	\begin{theorem}\label{t12}  
	Let $d=m=1$, let $k\in\N_0$ and assume that $a$ and $b$ satisfy 
		\begin{align}
				\text{There exist }  -\infty = \xi_0 < \xi_1 < \dots < \xi_{k} < \xi_{k+1} = \infty  \text{ such that} \hspace*{1.35cm}& \tag{A1}\\
				a \text{ is Lipschitz on } (\xi_{i-1}, \xi_i)  \text{ for }i=1,\dots,k+1,  \hspace*{3.65cm}&\notag\\[.1cm]
				b \text{ is Lipschitz on }\R \text{ and }b(\xi_i)\neq 0 \text{ for }i=1,\dots,k,\hspace*{3.15cm}& \tag{A2}\\
				a \text{ and } b \text{ have Lipschitz derivatives  on } (\xi_{i-1}, \xi_i) \text{ for }i=1,\dots,k+1.
			\hspace*{.18cm}&\tag{A3}	
	\end{align}
Then there exists a sequence of measurable functions $u_n\colon \R^n\to \R$, $n\in\N$, such that 
	for every $p\in [1,\infty)$,
		\[
			\EE \bigl[|X_1-u_n(W_{1/n},W_{2/n},\dots, W_1)|^p\bigr]^{1/p} \preceq n^{-3/4}.
		\]
	\end{theorem}

\begin{rem}\label{rem6}	
	Similar to the construction of a transformed Euler-Maruyama scheme, see Remark~\ref{rem5}(ii), the approximation  $u_n(W_{1/n},W_{2/n},\dots, W_1)$ in Theorem~\ref{t12} is obtained by transforming the original equation~\eqref{ee1} into a new SDE with coefficients that have sufficient regularity to obtain an $L_p$-error rate for the corresponding Milstein approximation and then backtransforming this  Milstein approximation to obtain an approximation for the original SDE.
	
	It is so far open whether in the setting of Theorem~\ref{t12} an $L_p$-error rate of at least $3/4$ is also achieved by the Milstein approximation $\mil_{n,1}$ of $X_1$ itself. See, however, Remark~\ref{rem8}.
\end{rem}

In contrast to the $L_p$-error rate $1$ for the Milstein scheme provided by Theorem~\ref{t2}, we only obtain an $L_p$-error rate of at least $3/4$ in Theorem~\ref{t12}. The following result from~\cite{ELL24,MGY2023} shows that this is neither a shortcoming of the approximation method nor of the associated technique of proof used in~\cite{MGY2022}. If the drift coefficient has at least one discontinuity then the rate $3/4$ can not be improved in general.

	\begin{theorem}\label{t13}  
	Let $d=m=1$, let $k\in\N$,  let $b=1$ and assume that $a$ satisfies 
		\begin{itemize}
		\item[(i)] there exist $-\infty=\xi_0<\xi_1<\ldots<\xi_k<\xi_{k+1}=\infty$ such that 
		$a$ has a bounded Lipschitz derivative on  $(\xi_{i-1}, \xi_i)$  for $i=1, \ldots, k+1$,
		\item[(ii)] there exists $i\in\{1, \ldots, k\}$ such that  $a(\xi_i+)\not=a(\xi_i-)$.
	\end{itemize}
Then 
	\[
	e_n^{(1)} \succeq n^{-3/4}.
	\]
\end{theorem}

Note that in the setting of Theorem~\ref{t13}, the assumptions (A1) to (A3) in Theorem~\ref{t12} are satisfied.
Combining Theorem~\ref{t12} and Theorem~\ref{t13} we thus obtain the following complexity result.

	\begin{theorem}\label{t14}  
		Let $d=m=1$, let $k\in\N$, let $b=1$ and assume that $a$ satisfies 
	\begin{itemize}
		\item[(i)] there exist $-\infty=\xi_0<\xi_1<\ldots<\xi_k<\xi_{k+1}=\infty$ such that 
		$a$ has a bounded Lipschitz derivative on  $(\xi_{i-1}, \xi_i)$  for $i=1, \ldots, k+1$,
		\item[(ii)] there exists $i\in\{1, \ldots, k\}$ such that  $a(\xi_i+)\not=a(\xi_i-)$. 
	\end{itemize}
	Then, for every $p\in [1,\infty)$,
	\[
	e_n^{(p)} \asymp n^{-3/4}.
	\]
\end{theorem}

\begin{rem}\label{rem7}	
In~\cite{Y2022}, a transformed Milstein approximation, see Remark~\ref{rem6}, is combined with an adaptive step-size control similar to the one used for the Euler-Maruyama scheme in~\cite{NSS2019}, see Remark~\ref{rem5}(ii). The resulting method achieves an $L_p$-error rate of at least $1$. Thus, in the case of a discontinuous drift coefficient, the  $L_p$-error rate of adaptive methods can be drastically better than the best possible $L_p$-error rate of non-adaptive methods. This is in stark contrast to the case of globally Lipschitz continuous coefficients, see Remark~\ref{rem3}(i).
\end{rem}

\subsection{Complexity results for scalar SDEs with a drift coefficient of fractional Sobolev regularity}\label{sobolev}
 
For quantifying the regularity of the drift coefficient we use the space 
	\[
	W^{s,p} = \biggl\{ a\colon\R\to\R \text{ measurable }\Bigl |\,  \int_{\R}\int_{\R} \frac{|a(x)-a(y)|^p}{|x-y|^{1+sp}}\,dx\, dy < \infty \biggr\}
	\]
of functions $a$ with fractional Sobolev regularity  $s \in (0,1)$ and integrability exponent $p\in [1,\infty)$.
	
The following theorem is a recent result from~\cite{DGL22} on the performance of the Euler-Maruyama approximation in case of an SDE~\eqref{ee1} with  a drift coefficient $a\in W^{s,p} $  and additive noise, i.e. $b=1$. Note that in this setting, the Euler-Maruyama approximation coincides with the Milstein approximation, i.e.
\begin{equation}\label{eulermil}
\forall  n\in\N\colon\,\, \eul_{n,1} =\mil_{n,1}.
\end{equation}

\begin{theorem}\label{t15}
	Let $d=m=1$, let $s\in (0,1)$ and $p\in [2,\infty)$, let $b=1$ and assume that  $a$ is bounded and $a\in W^{s,p}$. Then
		\[
		\forall \eps \in (0,\infty)\colon\,\,	\EE\bigl[ |X_1- \eul_{n,1}|^p] \bigr]^{1/p}  \preceq n^{-(1+s)/2+\eps}. 	
		\]
\end{theorem}

\begin{rem}\label{rem7b}
\noindent	
\begin{itemize}
	\item[(i)] Under the assumptions in Theorem~\ref{t15}, the equation~\eqref{ee1} has a unique strong solution. This follows from a more general result on existence of solutions of SDEs in~\cite{V80}.
	\item [(ii)] For $p=2$ and a slightly different setting for the drift coefficient $a$, the error rate $(1+s)/2-$ from Theorem~\ref{t15} was previously obtained in~\cite{NS2021}. To be more precise, if $b=1$ and  $a=a_1+a_2$, where  $a_1\colon\R\to\R$ is bounded and has bounded continuous derivatives up to order $2$ and $a_2\colon\R\to\R$  is bounded, Lebesgue-integrable and satisfies $a_2\in  W^{s,2}$  then the $L_2$-error rate $(1+s)/2-$  is in~\cite{NS2021} shown to be achieved by the equidistant Euler-Maruyama scheme if $s\le 1/2$ and by the Euler-Maruyama scheme based on the discretization $t_i= i^2/n^2$, $i=0,1,\dots,n$, if $s>1/2$. 
\end{itemize}
\end{rem}

A well-known example for the setting of Theorem~\ref{t15} is provided by H\"older continuous functions with compact support. If $a\colon\R\to\R$ is H\"older continuous of order $\gamma\in (0,1]$ and has compact support then 
\[
	a\in  \bigcap_{0<s<  \gamma } {W}^{s,2},
	\]
see, e.g.~\cite[Example 2.17]{Ern21} and \cite[Lemma 5.1]{DNPV2012}.	
Furthermore, piecewise H\"older regularity can be captured as well. If $a$ is piecewise H\"older continuous of order $\gamma\in (0,1]$ and has compact support then 
	\[
	a\in  \bigcap_{0<s<  \gamma \wedge 1/2} {W}^{s,2},
	\] 
see e.g. \cite[Section 3.1]{Sickel2021} and \cite[Theorem 4.6.4/2]{rs96}.
	 This
shows that the setting of piecewise regularity in Chapter~\ref{piecewise} and the setting of fractional Sobolev regularity have a certain overlap. See, in particular, the following remark.

\begin{rem}\label{rem8}
Assume that $b=1$ and that $a$ has compact support and satisfies the conditions  in Theorem~\ref{t12}, i.e. $a$ has piecewise bounded Lipschitz derivatives. Then an $L_2$-error rate of at least $3/4$ is obtained by a transformed Milstein scheme, see Theorem~\ref{t12} and Remark~\ref{rem6}. On the other hand,  $a$ is in particular piecewise Lipschitz continuous and we thus obtain  $a\in {W}^{s,2}$ for every $s\in (0,1/2)$. Hence, by~\eqref{eulermil} and Theorem~\ref{t15},
	\[
	\forall \eps \in (0,\infty)\colon\,\,	\EE\bigl[ |X_1- \mil_{n,1}|^2] \bigr]^{1/2}=\EE\bigl[ |X_1- \eul_{n,1}|^2] \bigr]^{1/2}  \preceq n^{-3/4+\eps}. 	
	\]
We thus see that, up to an arbitrary small $\varepsilon$, the Milstein approximation itself achieves an $L_2$-error rate of at least $3/4$ in this case. 

Moreover, if, additionally, $a$ has at least one discontinuity then the assumptions in Theorem~\ref{t13} are satisfied and the $L_2$-error rate $3/4$ is optimal.
\end{rem}

The following result from~\cite{EMGY2024} shows that, at least for $s\in (1/2,1)$ and $p=2$, the error rate $(1+s)/2-$ from Theorem~\ref{t15} can in general  essentially not be improved. For $s\in (1/2,1)$ define
\[
h_s\colon \R\to\R, \, x\mapsto \frac{ 1}{ \ln(e+|x|) \,(1+|x|)^{1/2+s}}
\]
and note that $h_s$ is Lebesgue integrable and symmetric. Let 
\begin{equation}\label{rr4}
a_s\colon\R\to\R,\quad x\mapsto 2\int_0^\infty cos(xz)\, h_s(z)\,dz
\end{equation}
denote the Fourier transform of $h_s$. 

\begin{theorem}\label{t16}
Let $d=m=1$ and let $s\in (1/2,1)$. The function $a_s$ given by~\eqref{rr4}
is bounded, Lebesgue-integrable and satisfies $a_s\in W^{s,2}$. For the equation~\eqref{ee1} with drift coefficient $a_s$ and diffusion coefficient $b=1$,
	\[
	e^{(2)}_n \succeq 	\frac{1}{n^{(1+s)/2}\ln (n+1)}.
	\]
\end{theorem}

\begin{rem}\label{rem10} A slight modification of the  function $a_s$ is already used in~\cite{Altm2021} to prove a lower $L_2$-error bound for approximating Brownian occupation time functionals $\int_0^1 a(W_t)\, dt$ with $a\in W^{s,2}$ based on equidistant evaluation of $W$. 
\end{rem}	

Combining Theorem~\ref{t15} and Theorem~\ref{t16} we get the following complexity result. 

\begin{theorem}\label{t17}
Let $d=m=1$ and let $b=1$. For every $s\in (1/2,1)$ there exists a bounded, Lebesgue integrable drift coefficient $a\in W^{s,2}$ such that 
\[
\forall\eps\in (0,\infty)\colon\,\, n^{-(1+s)/2-\eps}	\preceq 	e^{(2)}_n \preceq n^{-(1+s)/2+\eps}. 
\] 	
\end{theorem}

\begin{rem}\label{rem11}
	We briefly comment on a fundamental  tool to obtain lower $L_p$-error bounds for strong approximations of SDEs, which is used in particular in the proofs of Theorem~\ref{t13} and Theorem~\ref{t16}. In both cases, the proof relies on a coupling of noise approach, which in turn is a particular instance of a principle from {\sl Information Based Complexity} termed {\sl radius of information} or, more loosely, {\sl fooling algorithms}, see, e.g~\cite[Chapter 4]{NW08}. In the context of strong approximation of SDEs, the idea of this principle is as follows. Assume that the equation~\eqref{ee1} has a unique strong solution. Given $n\in\N$ and sites $\pi = (t_1,\dots, t_n)\in [0,1]^n$ for evaluating  the driving Brownian motion $W$, consider a further Brownian motion $ W^\pi $ such that $W$ and $W^\pi $ are coupled at the points $t_1,\dots,t_n$, i.e. 
	\begin{equation}\label{fool}
	W(t_i) = W^\pi(t_i)\,\text{ for }\, i=1,\dots,n.
\end{equation}
 Let $X^\pi$ denote the strong solution of the equation~\eqref{ee1} with driving Brownian motion $W^\pi$ in place of $W$. By a well known representation theorem for SDEs, there exists a measurable function $\psi\colon C([0,1];\R)\to\R$ such that a.s.
	\begin{equation}\label{fool0}
		X_1 = \psi(W)\text{ and }X^\pi_1 = \psi(W^\pi).
		\end{equation}
	For any approximation $u(W_{t_1},\dots,W_{t_n})$ of $X_1$ based on the evaluation of $W$ at the points $t_1,\dots,t_n$ we then have, due to~\eqref{fool},
	\begin{equation}\label{fool1}
	u(W_{t_1},\dots,W_{t_n}) = u(W^\pi_{t_1},\dots,W^\pi_{t_n}),
	\end{equation}
i.e., the algorithm determined by the function $u$ can not distinguish between strong approximation of $X_1$ and strong approximation of $X^\pi_1$ on the basis of the information provided by evaluating the respective driving Brownian motion at the points $t_1,\dots,t_n$. Using~\eqref{fool0}  we obtain that there exists 
$\varphi\colon C([0,1];\R)\to\R$ such that
\begin{equation}\label{fool2}
	X_1 -u(W_{t_1},\dots,W_{t_n}) = \varphi (W)\text{ and } X^\pi_1-u(W^\pi_{t_1},\dots,W^\pi_{t_n})=\varphi(W^\pi).
\end{equation}	
By~\eqref{fool1},  ~\eqref{fool2} and  the triangle inequality we conclude that for every $p\in [1,\infty)$, 
		\begin{equation}\label{fool3}
			\begin{aligned}
	&	\EE\bigl[|X_1-u(W_{t_1},\dots,W_{t_n})|^p\bigr]^{1/p} \\
	&\qquad \qquad = \frac{1}{2} \Bigl( \EE\bigl[|X_1-u(W_{t_1},\dots,W_{t_n})|^p\bigr]^{1/p} + \EE\bigl[|X^\pi_1-u(W^\pi_{t_1},\dots,W^\pi_{t_n})|^p\bigr]^{1/p} \Bigr)\\
	& \qquad\qquad = \frac{1}{2} \Bigl( \EE\bigl[|X_1-u(W_{t_1},\dots,W_{t_n})|^p\bigr]^{1/p} + \EE\bigl[|X^\pi_1-u(W_{t_1},\dots,W_{t_n})|^p\bigr]^{1/p} \Bigr)\\
	& \qquad\qquad \ge \frac{1}{2} \EE\bigl[|X_1-X^\pi_1\bigr|^p]^{1/p}.
		\end{aligned}
	\end{equation}
In this way, the $L_p$-error of any approximation of $X_1$ based on evaluating $W$ at the points from $\pi$ is, up to a constant, bounded from below by the $L_p$-distance of $X_1$ and $X^\pi_1$. Clearly, $W^\pi$ should be chosen such  that the latter quantity is large. The particular approach used in the proof of Theorem~\ref{t13} and Theorem~\ref{t16} is to construct $W^\pi$ such that, conditioned on $W_{t_1},\dots,W_{t_n}$, the two processes $W$ and $W^\pi$ are independent. In the case $0=t_0<t_1<\dots <t_n=1$ this can be achieved by taking
\[
W^\pi = \overline W^\pi + U^\pi,
\]
where $\overline W^\pi$ denotes the piecewise linear interpolation of $W$ at the points $t_0,\dots,t_n$, $U^\pi$ coincides with the piecewise Brownian bridge process  $W- \overline W^\pi$ in distribution and $U^\pi$ and $W$ are independent.
	\end{rem}

\bibliographystyle{acm}
\bibliography{bibfile}

	\end{document}